\documentclass[12pt]{article}
\usepackage{amsfonts,amsbsy,amssymb,amsmath,amsthm}
\usepackage[toc,title,page]{appendix}
\usepackage{a4}
\usepackage{graphicx}
\usepackage[english]{babel}
\usepackage{datetime}
\usepackage[colorinlistoftodos]{todonotes}
\usepackage[numbers,sort]{natbib}

\newcommand{\ZZ}{\mathbb{Z}}
\newcommand{\CC}{\mathbb{C}}
\newcommand{\RR}{\mathbb{R}}
\newcommand{\NN}{\mathbb{N}}
\newcommand{\ab}{{\boldsymbol{a}}}
\newcommand{\bb}{{\boldsymbol{b}}}
\newcommand{\cb}{{\boldsymbol{c}}}
\newcommand{\vb}{{\boldsymbol{v}}}
\newcommand{\db}{{\boldsymbol{d}}}
\newcommand{\pb}{{\boldsymbol{p}}}
\newcommand{\Ab}{{\boldsymbol{A}}}
\newcommand{\Bb}{{\boldsymbol{B}}}
\newcommand{\Db}{{\boldsymbol{D}}}
\newcommand{\Fb}{{\boldsymbol{F}}}
\newcommand{\Ib}{{\boldsymbol{I}}}
\newcommand{\Zerob}{{\boldsymbol{0}}}

\newcommand{\deltab}{{\boldsymbol{\delta}}}
\newcommand{\supp}{\operatorname{supp}}

\newcommand{\cT}{{\cal T}}
\newcommand{\cL}{\mathcal{L}}
\newcommand{\cS}{\mathcal{S}}
\newcommand{\diag}{\operatorname{diag}}
\newcommand{\id}{\operatorname{id}}
\newcommand{\bPhi}{\boldsymbol{\Phi}}
\newcommand{\dimV}{m}
\newcommand{\itVar}{n}
\newcommand{\vfield}{\mathcal{V}}

\usepackage{cleveref}

\newtheorem{Theorem}{Theorem}
\newtheorem{Definition}[Theorem]{Definition}
\newtheorem{Proposition}[Theorem]{Proposition}
\newtheorem{Lemma}[Theorem]{Lemma}

\newtheorem{Remark}[Theorem]{Remark}
\title{Hermite multiwavelets for manifold-valued data}
\author{
  {Mariantonia Cotronei}\thanks{DIIES, Universit\`a Mediterranea di
    Reggio Calabria, Via Graziella loc. Feo di Vito, 89122 Reggio
    Calabria, Italy. \texttt{mariantonia.cotronei@unirc.it}}  
\and 
{Caroline Moosm\"uller}\thanks{Department of Mathematics, University of California, San Diego, 9500 Gilman Drive, La Jolla, CA 92093, USA. \texttt{cmoosmueller@ucsd.edu}} 
\and 
{Tomas Sauer}\thanks{Lehrstuhl f\"ur Mathematik mit Schwerpunkt Digitale
  Signalverarbeitung \& FORWISS, Universit\"at Passau,
  Fraunhofer IIS Research Group on Knowledge Based Image
  Processing, Innstr.~43,
  94032 Passau, Germany. \texttt{tomas.sauer@uni-passau.de}}
\and
{Nada Sissouno}\thanks{Department of Mathematics, Technical University
  of Munich, Boltzmannstra\ss e 3, 85748 Garching,
  Germany. \texttt{sissouno@ma.tum.de}}
}

\date{}

\begin{document}
%%%%%%%%%%%%%%%%%%%%%%%%%%%%%%%%%%%%%%%%
\maketitle

%%%%%%%%%% ABSTRACT %%%%%%%%%%%%%%%%%%%%%%%%%%%%%%%

\begin{abstract}
In this paper we present a construction of interpolatory Hermite multiwavelets for functions that take values in nonlinear geometries such as Riemannian manifolds or Lie groups. We rely on the strong connection between wavelets and subdivision schemes to define a prediction-correction approach based on Hermite subdivision schemes that operate on manifold-valued data. The main result concerns the decay of the wavelet coefficients: We show that our manifold-valued construction essentially admits the same coefficient decay as linear Hermite wavelets, which also generalizes results on manifold-valued scalar wavelets.

\end{abstract}

%65T60; Numerical methods for wavelets
%65D15; Algorithms for approximation of functions
% 41A25: Rate of convergence, degree of approximation
% 53A99 Differential geometry

%%%%%%%%%%%%%%%%%%%%%%%%%%%%%%%%%%%%%%%%%%%%%%%%%%%%
%%%%%%%%%%%%%%%%%%%%%%%%%%%%%%%%%%%%%%%%%%%%%%%%%%%
%%%%%%% INTRO %%%%%%%%%%%%%%%%%%%%%%%%%%%%%%
\section{Introduction}

%% waevelet transforms
Wavelets are one of the most important tools for the analysis of signals and images, as they allow to study local properties of functions at different resolutions.
In the last decades a lot of different types of one or multidimensional wavelets as well as their properties have been studied; see \cite{daubechies92,mallat09:_wavel_tour_signal_proces} for an overview.
The literature on wavelet transforms for functions that take values in nonlinear geometries, such as Riemannian manifolds or Lie groups, is not as exhaustive as in the linear case. In the manifold setting, the aim is to construct processes which are intrinsic to the underlying geometry, for example by preserving invariances with respect to certain transformation groups.

The idea of considering geometric data goes back to \cite{Rahman05}, and has led to a series of results concerning convergence and smoothness of subdivision schemes, starting with the work of \cite{wallner05,wallner06}, the coefficient decay for interpolatory wavelets \cite{grohs09}, and the definability and stability of multiscale transforms \cite{grohs12,grohs10b}.

In this paper, we aim at extending this line of research by defining and analyzing \emph{multiwavelets} for manifold-valued data. Linear multiwavelets are a generalization of classical (scalar) wavelets and are obtained by allowing several function in the construction of multiresolution analyses. They are based on a \emph{multi-scaling function} that satisfies a \emph{vector refinement equation} with matrix-valued rather than scalar coefficients. 
Multiwavelets can have advantageous properties, for example, for constructing bases with short support and high approximation \cite{keinert04}.

This paper focuses on multiwavelets of Hermite-type, meaning that the multi-scaling functions satisfy Hermite conditions \cite{cotronei17,cotronei19,strela1995}. 
Such wavelet systems can find applications in contexts  where Hermite data  need to be processed, for example for  compression or denoising reasons.

In particular, starting from an \emph{interpolatory Hermite subdivision scheme} reproducing elements in a given space, for example the space of polynomials or exponential functions, it is  always possible to realize a {\em biorthogonal wavelet system}, where the associated  wavelet operator  possesses the property of ``cancelling" those elements \cite{cotronei17}. This is the usually required vanishing moment property  assuring good compression capabilities to the wavelet system.

We use the mentioned tight connection between subdivision schemes and wavelets to obtain manifold-valued Hermite wavelet schemes, using the construction presented in \cite{moosmueller16,moosmueller17}. This works in a similar fashion as the scalar constructions of \cite{grohs09,grohs12}.

The main result of this paper is a wavelet coefficient decay property of such manifold-valued wavelets, which mimics the linear case \cite{cotronei19} and can be considered as an extension of \cite{grohs09} to Hermite-type interpolatory wavelets. 

The paper is organized as follows. In \Cref{sec:preliminaries} and \Cref{sec:linSubDiv} we introduce the linear tools necessary to construct Hermite-type wavelets, mainly focusing on Hermite subdivision schemes. \Cref{sec:results} introduces linear Hermite multiwavelets based on \cite{cotronei17}. We reinterpret their constructions in terms of operators rather than symbols, highlighting the similarities with the scalar multiscale transforms of \cite{grohs12}. \Cref{sec:manifoldSubDiv} introduces our Hermite prediction-correction scheme for manifold-valued data, which is a direct generalization of \cite{cotronei17} and makes use of natural tools in nonlinear geometries such as the exponential map and the parallel transport operator. 
In this section we also prove that the wavelet coefficients at level $n$ decay as $2^{-2n}$ for dense enough input data, showing that manifold-valued Hermite wavelets have similar properties as their linear counterparts \cite{cotronei19}.

%%%%%%%%%%%%%%%%%%%%%%%%%%%%%%%%%%%%%%%%%%%%%%%%%%
% PRELIMINARIES
%%%%%%%%%%%%%%%%%%%%%%%%%%%%%%%%%%%%%%%%%%%%%%%%%%

%%%%%%%%%%%%%%%%%%%%%%%%%%%%%%%%%%%%%%%%%%
\section{Preliminaries}\label{sec:preliminaries}
%%%%%%%%%%%%%%%%%%%%%%%%%%%%%%%%%%%%%%%%%%
%
% short intro
In this paper we are concerned with wavelets for functions $f:\RR \to M$, where $M$ is a manifold. The main examples of manifolds we consider are surface in $\RR^{\dimV}$ and Lie groups. To construct wavelets for manifold-valued functions, we also include information about the first derivatives $f'$.

In the linear version of this problem, the data are of the
form $(f(x),f'(x))^T \in \RR^{\dimV}\times \RR^{\dimV}$ for $x\in \RR$. To simplify notation, we denote by $V=\RR^{\dimV}$, so that the data lies in $V^2$. Throughout this text, $\dimV$ always denotes the dimension of $V$.

%% vector/matrix notation
Elements in $V^2$ are denoted by bold lower case letters $\pb$. We are also concerned with $L(V)^{2\times 2}$, where $L(V)$ is the space of all linear functions $V\to V$. Elements of $L(V)^{2\times 2}$ are denoted by bold upper case letters $\Ab$.
%% sequence spaces
The space of all vector-valued sequences $\ZZ \to V^2$ is denoted by $\ell(\ZZ,V^2)$. Elements of $\ell(\ZZ,V^2)$ are again denoted by bold lower case letters $\pb = (\pb_j: j\in \ZZ)$. 
We also consider the matrix-valued sequence space $\ell(\ZZ,L(V)^{2\times 2})$. Elements of this space are again denoted by bold upper case letters $\Ab = (\Ab_j: j\in \ZZ)$.

%% norms
We introduce norms on $\ell(\ZZ,V^2)$ and $\ell(\ZZ,L(V)^{2\times 2})$:
\begin{align}\label{eq:norm}
    \|\pb\|_{\infty} &= \sup_{j\in\ZZ}|\pb_j|_{\infty} \\ \nonumber
    \|\Ab\|_{\infty} &= \sup_{j\in\ZZ}|\Ab_j|_{\infty},
\end{align}
where $|\cdot|_{\infty}$ denotes the infinity-norm in $V^2$ (resp.\ $L(V)^{2\times 2}$). The space consisting only of bounded sequences with respect to the norms (\ref{eq:norm}) are denoted by $\ell_{\infty}(\ZZ,V^2)$ and $\ell_{\infty}(\ZZ,L(V)^{2\times 2})$.
We further consider $\ell_{0}(\ZZ,L(V)^{2\times 2})$, which is the space of finitely supported sequences in $\ell(\ZZ,L(V)^{2\times 2})$.

%% function spaces
By $C(\RR,V)$ we denote the space of continuous functions $\RR \to V$, while $C_u(\RR,V)$ denotes the space of uniformly continuous and bounded functions. We further consider the space of continuously differentiable functions $C^1(\RR,V)$ and the space $C^1_u(\RR,V)$ of functions $f\in C^1(\RR,V)$ with $f'\in C_u(\RR,V)$.

%% Decomposition and reconstruction
The decomposition and
reconstruction of data using filter banks is closely related
to wavelets and subdivision schemes. A detailed discussion of the connection of filter banks and wavelets, especially in the setting of biorthogonal wavelets that we analyze, can be found in \cite{vetterli92,vetterli95}.

We consider \emph{filters} or \emph{masks} $\Ab^{[n]}\in \ell_0(\ZZ,L(V)^{2\times 2})$, $\itVar \in \NN$, of the form
\begin{equation}\label{eq:typeMasks}
    \Ab^{[\itVar]} = 
    \left(
    \begin{array}{cc}
        \ab_{00}^{[\itVar]} & \ab_{01}^{[\itVar]} \\[0.4em]
        \ab_{10}^{[\itVar]}  &  \ab_{11}^{[\itVar]}
    \end{array}
    \right),
\end{equation}
where $\ab_{00}^{[\itVar]},\ab_{10}^{[\itVar]},\ab_{01}^{[\itVar]},\ab_{11}^{[\itVar]} \in \ell(\ZZ,\RR)$. The entries of $\Ab^{[\itVar]}$ in \cref{eq:typeMasks} are to be understood as $\ab_{00}^{[\itVar]}\cdot \Ib$, etc., where $\Ib$ denotes the identity matrix. Through this form of $\Ab^{[\itVar]}$, results for Hermite subdivision schemes with $V=\RR$ can be directly applied to our set up.
An important mask is the \emph{delta sequence} $\deltab=(\deltab_j:j\in\ZZ)$ given by 
$\deltab_0=\Ib$ and $\deltab_j=\Zerob$ for $j\in\ZZ\backslash \{0\}$.

Let $\pb \in \ell(\ZZ,V^2)$ and $j\in \ZZ$. Given a mask $\Ab^{[n]}$,
the associated \emph{reconstruction} or \emph{subdivision operator of level $\itVar$}, $\mathcal{S}_{\Ab^{[\itVar]}}:\ell(\ZZ,V^2) \to \ell(\ZZ,V^2)$, is given by
\begin{equation}
\label{eq:subdivisionOperator}
    (\mathcal{S}_{\Ab^{[\itVar]}}\pb)_j = \sum_{k\in \ZZ} \Ab^{[\itVar]}_{j-2k}\pb_k,
\end{equation}
while the \emph{decomposition} or \emph{wavelet operator} $\mathcal{D}_{\Ab^{[n]}}$ is given by
$$(\mathcal{D}_{\Ab^{[n]}}\pb)_j=\sum_{i\in \mathbb{Z}}\Ab_{i-2j}^{[n]}\pb_i.$$
We also need the \emph{shift operator} $\mathcal{L}:\ell(\ZZ,V^2) \to \ell(\ZZ,V^2)$ defined as
\begin{equation}\label{eq:shift}
    (\mathcal{L}\pb)_i=\pb_{i+1}.
\end{equation}
The reconstruction and decomposition operators satisfy the following well-known properties:
\begin{equation}
\label{eq:sub-wavelet-shift-commute}    
\mathcal{S}_{\Ab^{[\itVar]}}\mathcal{L}
=
\mathcal{L}^2\mathcal{S}_{\Ab^{[\itVar]}}
\quad \text{and}\quad
\mathcal{D}_{\Ab^{[\itVar]}}\mathcal{L}^2
=
\mathcal{L}\mathcal{D}_{\Ab^{[\itVar]}}
\end{equation}

%%%%%%%%%%%%%%%%%%%%%%%%%%%%%%%%%%%%%%%%%%%%%%%%%%%%%%%%%%%%
% HERMITE SUBDIVISION SCHEMES
%%%%%%%%%%%%%%%%%%%%%%%%%%%%%%%%%%%%%%%%%%%%%%%%%%%%%%%%%%%%
%%%%%%%%%%%%%%%%%%%%%%%%%%%%%%%%%%%%%%%%%%%%%%
\section{Linear Hermite subdivision schemes}\label{sec:linSubDiv}
%%%%%%%%%%%%%%%%%%%%%%%%%%%%%%%%%%%%%%%%%%%%%%
%
%
%% subdivision scheme
%
Consider a sequence of finitely supported masks $(\Ab^{[\itVar]}:\itVar \geq 0)$. A \emph{linear Hermite subdivision scheme} $S(\Ab^{[\itVar]}:\itVar \geq 0)$ is the iterative procedure of constructing sequences $\pb^{[\itVar]}$ from an initial sequence $\pb^{[0]}$ via the rule
\begin{equation}
\label{eq:subdivisionScheme}
\Db^{\itVar+1}\pb^{[\itVar+1]} = \mathcal{S}_{\Ab^{[\itVar]}}\Db^{\itVar} \pb^{[\itVar]},
\quad \itVar \in \NN.
\end{equation}
%
% matrix D
%
Here $\Db$ denotes the $\Db = \diag\left(1,1/2\right) \in L(V)^{2\times 2}$, where a constant $C$ is to be understood as $C\cdot \Ib$.
Since we associate $\pb^{[\itVar]}$ with pairs of function and derivative each evaluated on the grid $2^{-\itVar}\ZZ$, the matrix $\Db$ and its powers arise because of the chain rule.

%% stationary
Schemes of the form \cref{eq:subdivisionScheme} are often called \emph{level-dependent} as opposed to \emph{stationary}. In stationary subdivision $\Ab^{[\itVar]}=\Ab$ is satisfied for a fixed mask $\Ab$, i.e.\ the mask does not depend on the iteration level $\itVar$.

%% interpolatory schemes
In this paper we are mostly concerned with \emph{interpolatory schemes}: A scheme satisfying \cref{eq:subdivisionScheme} is called interpolatory if $\pb^{[\itVar+1]}_{2j}=\pb^{[\itVar]}_j$ for $j\in \ZZ, \itVar \in \NN$. This property relates to the sequence of masks $(\Ab^{[\itVar]},\itVar \in \NN)$ satisfying $\Ab^{[\itVar]}_{2j}=\Db \deltab_j, j\in \ZZ, n\in \NN$.
%where $\deltab_j$ are the elements of the delta sequence $\deltab$.
In terms of operators the interpolation property can be written as 
$\mathcal{D}_{\deltab}\mathcal{S}_{\Ab^{[n]}}=\Db$.

%% convergence
A Hermite subdivision scheme is called $C^1$-convergent if for every initial data $\pb^{[0]} \in \ell_{\infty}(\ZZ,V^2)$ there exists a function $\bPhi = [\Phi_k]_{k=0}^1:\RR \to V^2$ such that the sequence $\pb^{[\itVar]}$ satisfies
\begin{equation*}
    \lim_{n\to \infty}\sup_{j\in\ZZ}
    \left| \pb^{[\itVar]}_j - \bPhi\left(\frac{j}{2^{\itVar}}\right)\right|_{\infty} = 0,
\end{equation*}
and where $\Phi_0 \in C^1_u(\RR,V)$ with $\Phi_0' = \Phi_1$. We further assume that there exists at least one sequence $\pb^{[0]} \in \ell(\ZZ,V^2)$ such that the resulting limiting function satisfies $\bPhi\neq 0$.
%% Convergence results
Results on the convergence of linear Hermite subdivision schemes can be found, for example, in \cite{dubuc05,merrien12,dubuc06,dyn95,han05,han03b} for the stationary case, in \cite{conti17a,jeong19} in the level-dependent case, and in \cite{conti14,moosmueller20,moosmueller21,han20} for smoothness of high order.

% Basic limit function
When applying a $C^1$-convergent scheme to the delta sequence as initial data it converges to the so-called \emph{basic limit function}
\begin{equation*}
    \Fb =\left(
    \begin{array}{cc}
        \Phi_0 & \Phi_1\\[0.4em]
        \Phi'_0  &  \Phi'_1
    \end{array}
    \right),
\end{equation*}
see \cite{dubuc05} for the case of Hermite schemes.
If we consider $C^1$-convergent schemes starting at level $\ell$, i.e. $S(\Ab^{[\itVar+\ell]}:\itVar \geq 0)$ for $\ell\geq 0$  applied to the delta sequence, we obtain a sequence of 
basic limit functions $\Fb^{[\ell]}$ with $\Fb^{[0]}=\Fb$.
The basic limit functions at different levels are connected via a refinement equation, which allows to use them for the construction of multiresolution analyses \cite{cotronei17,cotronei19}.

%Data and Sampling
Closely related to the convergence of subdivision schemes and the refinement property is the property of
reproducing certain spaces \cite{conti16,conti17a,jeong17,merrien12}.
Here we consider Hermite subdivision schemes that reproduce at least a $2$-dimensional space of polynomials and/or exponentials. Since reproduction of constants is a necessary condition for convergence, the space to be reproduced should either contain
\begin{equation}\label{eq:space_W}
\operatorname{span}\{ 1,x\} \quad \text{or} \quad \operatorname{span}\{1,e^{\lambda x}\},
\end{equation}
where $\lambda \in
\CC \setminus \{ 0 \}$. Some examples of Hermite schemes reproducing such spaces can be found in \cite{conti16,conti15,conti14,jeong17,jeong19}.
In the following, we write $W$ to mean either one of the spaces in \eqref{eq:space_W}.

The reproduction property can be formulated in
terms of the \emph{spectral condition} \cite{conti16,dubuc09,merrien12} or \emph{sum rules} \cite{han05}:
\begin{equation*}\label{eq:spectralcond}
    \mathcal{S}_{\Ab^{[n]}} \Db^{n} \vb^{[n]}_f=\Db^{n+1} \vb^{[n+1]}_f, \qquad f
    \in W, \, n\in \NN.
  \end{equation*}
where $\vb^{[n]}_f$ is the vector-valued sequences
\begin{equation*}%\label{eq:vfseq}
 \vb^{[n]}_{f;j}= 
 \left(\begin{array}{c}
        f(2^{-n}j)\\
        f'(2^{-n}j)
        \end{array}
\right),
\qquad j \in \ZZ.
\end{equation*}
defined by a function $f\in C^1(\RR)$.

%%%%%%%%%%%%%%%%%%%%%%%%%%%%%%%%%%%%%%%%%%%%%%%%%%%%%%%%%%%%%%%%%%%
\section{Linear wavelets from interpolatory Hermite subdivision schemes}
\label{sec:results}
%%%%%%%%%%%%%%%%%%%%%%%%%%%%%%%%%%%%%%%%%%%%%%%%%%%%%%%%%%%%%%%%%%%
In \cite{cotronei17} multiwavelets are constructed from linear Hermite subdivision schemes and \cite{cotronei19} provides an estimate on the wavelet coefficient decay. These papers rely on the \emph{symbol} of the matrix mask, i.e., the matrix-valued Laurent polynomials
\begin{equation}\label{eq:symbol}
 \Ab^{[n]}(z):=\sum_{k\in\ZZ}\Ab_k^{[n]}z^k, \quad   z\in\CC,
\end{equation}
and $\Ab^{[n]}\in\ell_0(\ZZ,L(V)^{2\times 2})$.
To generalize the results of \cite{cotronei17,cotronei19} to the manifold-valued case we rewrite the necessary constructions in terms of operators rather than symbols.

We consider sets of level-dependent filters $\{\Ab^{[n]},\Bb^{[n]},\widetilde\Ab^{[n]},\widetilde \Bb^{[n]}, \, n\in \NN\}$, where $\widetilde\Ab^{[n]}$ and $\widetilde \Bb^{[n]}$ are the filters associated to the decomposition of data and $\Ab^{[n]}$ as well as $\Bb^{[n]}$ denote filters associated to the reconstruction.

%%%%%%%%%%%%%%%%%%%%%%%%%%%%%%%%%%%%%%%%%%
% define biorthogonal system
%%%%%%%%%%%%%%%%%%%%%%%%%%%%%%%%%%%%%%%%%%
\begin{Definition}\label{def:biorthogonal}
Given a set of level-dependent filters $\{\Ab^{[n]},\Bb^{[n]},\widetilde\Ab^{[n]},\widetilde \Bb^{[n]}, \, n\in \NN\}$ we say that they form a \emph{biorthogonal system} if the following conditions are satisfied: 
\begin{eqnarray*}
    &&\mathcal{D}_{(\mathbf{\widetilde{A}}^{[n]})^T}\mathcal{S}_{\Ab^{[n]}}=
    \mathcal{D}_{(\mathbf{\widetilde{B}}^{[n]})^T}\mathcal{S}_{\Bb^{[n]}}=
    \operatorname{id},\\
&&\mathcal{D}_{(\mathbf{\widetilde{A}}^{[n]})^T}\mathcal{S}_{\Bb^{[n]}}=
\mathcal{D}_{(\mathbf{\widetilde{B}}^{[n]})^T}\mathcal{S}_{\Ab^{[n]}}=0,
\end{eqnarray*}
for all $n\in \NN$.
\end{Definition}
%%%%%%%%%%%%%%%%%%%%%%%%%%%%%%%%%%%%%%%%%

%\begin{Remark}
The biorthogonal system conditions of \Cref{def:biorthogonal} are exactly the biorthogonal system conditions formulated in terms of symbols in \cite{cotronei17}, as proved in the following Proposition. %\Cref{thm:biorthogonal_symbol}.
%\end{Remark}
\begin{Proposition}
\label{thm:biorthogonal_symbol}
The biorthogonal system conditions of \Cref{def:biorthogonal} are exactly the biorthogonal system conditions formulated in terms of symbols in \cite{cotronei17}.
\end{Proposition}
%%%%%%%%%%%%%%%%%%
\begin{proof}
The biorthogonal system conditions in terms of symbols of 
\cite[Eq.\ (6)]{cotronei17} are:
\begin{align*}
&({\widetilde{\Ab}}^{[n]})^{\sharp}(z)\Ab^{[n]}(z)+({\widetilde{\Ab}}^{[n]})^{\sharp}(-z)\Ab^{[n]}(-z)=2I, \\
&({\widetilde{\Ab}}^{[n]})^{\sharp}(z)\Bb^{[n]}(z)+({\widetilde{\Ab}}^{[n]})^{\sharp}(-z)\Bb^{[n]}(-z)=0, \\
&({\widetilde{\Bb}}^{[n]})^{\sharp}(z)\Ab^{[n]}(z)+({\widetilde{\Bb}}^{[n]})^{\sharp}(-z)\Ab^{[n]}(-z)=0, \\
&({\widetilde{\Bb}}^{[n]})^{\sharp}(z)\Bb^{[n]}(z)+({\widetilde{\Bb}}^{[n]})^{\sharp}(-z)\Bb^{[n]}(-z)=2I.
\end{align*}
where $({\widetilde{\Ab}}^{[n]})^{\sharp}(z)=(\Ab^{[n]})^T(z^{-1})$ (see the definition of the symbol \eqref{eq:symbol}).
We show that the first condition is the same as our first operator condition (\Cref{def:biorthogonal}); the rest can be proved analogously.

We compute the symbol from the first equation:
\begin{align*}
2I& =\sum_{i,j}(\widetilde{\Ab}^{[n]})^T_iz^{-i}\Ab^{[n]}_jz^{j}
	+\sum_{i,j}(\widetilde\Ab^{[n]})^T_i(-z)^{-i}\Ab^{[n]}_j(-z)^{j}\\
& =\sum_{i,j}(1+(-1)^{i+j})(\widetilde\Ab^{[n]})^T_i\Ab^{[n]}_jz^{j-i}\\
&=\sum_k\Big(\sum_i(1+(-1)^k)(\widetilde\Ab^{[n]})^T_i\Ab^{[n]}_{i+k}\Big)z^k
\end{align*}
This implies
\begin{equation*}
\sum_i(1+(-1)^k)(\widetilde\Ab^{[n]})^T_i\Ab^{[n]}_{i+k}
	=\begin{cases}
	2I & \text{if } k=0,\\
	0 & \text{if } k\neq 0.
	\end{cases}
\end{equation*}
In particular
\begin{equation}\label{case}
\sum_i(\widetilde\Ab^{[n]})^T_i\Ab^{[n]}_{i+k}
	=\begin{cases}
	I & \text{if } k=0,\\
	0 & \text{if } k\neq 0 \text{ and } k \text{ is even.}
	\end{cases}
\end{equation}
Now the equation using operators is
\begin{align*}
(\mathcal{D}_{(\mathbf{\widetilde{A}}^{[n]})^T}\mathcal{S}_{\Ab^{[n]}}\cb)_j &=
\sum_{i\in \mathbb{Z}}(\mathbf{\widetilde{A}}^{[n]})^T_{i-2j}(\mathcal{S}_{\Ab^{[n]}}\cb)_i
=\sum_{i,k}(\mathbf{\widetilde{A}}^{[n]})^T_{i-2j}\Ab^{[n]}_{i-2k}\cb_k\\
&=\sum_k\Big(\sum_i(\mathbf{\widetilde{A}}^{[n]})^T_{i-2j}\Ab^{[n]}_{i-2k} \Big)\cb_k
=\sum_k\Big(\sum_r(\mathbf{\widetilde{A}}^{[n]})^T_{r}\Ab^{[n]}_{r+2(j-k)} \Big)\cb_k
\end{align*}
Applying \eqref{case} we see that
\begin{equation*}
(\mathcal{D}_{(\mathbf{\widetilde{A}}^{[n]})^T}\mathcal{S}_{\Ab^{[n]}}\cb)_j=\cb_j.
\end{equation*}
Thus $\mathcal{D}_{(\mathbf{\widetilde{A}}^{[n]})^T}\mathcal{S}_{\Ab^{[n]}} = \id$.
\end{proof}
%%%%%%%%%%%%%%%%%%%%%%%%

\begin{Remark}
From the biorthogonal filter conditions (\Cref{def:biorthogonal}), it follows that if ${\Ab^{[n]}}$ satisfies the $W$-spectral condition then
$\widetilde{{\Bb}}^{[n]}$ satisfies the $W$-vanishing moment condition, i.e. elements of $W$ are canceled in the decomposition of data: $\mathcal{D}_{(\widetilde{{\Bb}}^{[n]})^T}\Db^{n+1}\vb^{[n+1]}_{f;k}=0$ for $f\in W$. 
 
Indeed, if  $ f\in W$, then
\begin{equation*}
 \mathcal{S}_{\Ab^{[n]}}\Db^n\vb^{[n]}_{f;k}=\Db^{n+1}\vb^{[n+1]}_{f;k} \implies \mathcal{D}_{(\mathbf{\widetilde B}^{[n]})^T}\Db^{n+1}\vb^{[n+1]}_{f;k}=\mathcal{D}_{(\mathbf{\widetilde B}^{[n]})^T}\mathcal{S}_{\Ab^{[n]}}\Db^{n}\vb^{[n]}_{f;k}=0.
 \end{equation*}
 
See \cite{cotronei17} for more details on the relation between spectral and vanishing moment conditions. 
 \end{Remark}

For a given level-dependent biorthogonal wavelet system $\{\Ab^{[n]},\Bb^{[n]},\widetilde\Ab^{[n]},\widetilde \Bb^{[n]}, \, n\in \NN\}$ we rewrite the discrete wavelet transform formula, for the decomposition and the reconstruction, in terms of the respective operators:
%%%%%%%%%%%%%%%%%%%%%%
% define decomposition
%%%%%%%%%%%%%%%%%%%%%%%
\begin{Definition}\label{def:dec}
Let $N\in \NN$ and $\cb^{[N]} \in \ell(\ZZ, V^2)$. For $n=N-1,\ldots,0$,  the \emph{decomposition scheme}  reads as
\begin{align*}
\cb^{[n]}&=\mathcal{D}_{(\mathbf{\widetilde{A}}^{[n]})^T}\cb^{[n+1]},\\
\db^{[n]}&=\mathcal{D}_{(\mathbf{\widetilde{B}}^{[n]})^T}\cb^{[n+1]}.
\end{align*}
\end{Definition}
%%%%%%%%%%%%%%%%%%%%
Repeated application of the decomposition scheme leads to coarse data $\cb^{[0]}$ and wavelet coefficients $\db^{[0]},\ldots,\db^{[N-1]}$. One can reconstruct the data $\cb^{[n]}$ via the reconstruction scheme:
%%%%%%%%%%%%%%%%%%%%
% define reconstruction
%%%%%%%%%%%%%%%%%%%%
\begin{Definition}
Let $N\in \NN$ and $\cb^{[0]}, \db^{[0]}, \ldots, \db^{[N-1]}  \in \ell(\ZZ, V^2)$. For $n=0, \ldots,N-1$, the 
 \emph{reconstruction scheme} reads as:
\begin{equation*}
\cb^{[n+1]}=\mathcal{S}_{\Ab^{[n]}}\cb^{[n]}+\mathcal{S}_{\Bb^{[n]}}\db^{[n]}.
\end{equation*}
\end{Definition}

The reconstruction of $\cb^{[n]}$, $n=1,\ldots, N$, is called
\emph{perfect reconstruction}, if 
\begin{equation}\label{eq:perfect_reconstr}
\mathcal{S}_{\Ab^{[n]}}\mathcal{D}_{(\mathbf{\widetilde{A}}^{[n]})^T}+\mathcal{S}_{\Bb^{[n]}}\mathcal{D}_{(\mathbf{\widetilde{B}}^{[n]})^T}
=\operatorname{id}.
\end{equation}
for all $n$.

Using the biorthogonality conditions (\Cref{def:biorthogonal}),  we may write the decomposition scheme in the following way (compare \cite[p.3, eq. (5)]{grohs12}):
\begin{align}\nonumber
\cb^{[n]}&=\mathcal{D}_{(\mathbf{\widetilde{A}}^{[n]})^T}\cb^{[n+1]},\\\nonumber
\db^{[n]}&=\mathcal{D}_{(\mathbf{\widetilde{B}}^{[n]})^T}\Big(\cb^{[n+1]}-\mathcal{S}_{\Ab^{[n]}}\cb^{[n]}\Big)\\ \label{dec2}
&=\mathcal{D}_{(\mathbf{\widetilde{B}}^{[n]})^T}
\Big(\operatorname{id}-\mathcal{S}_{\Ab^{[n]}}\mathcal{D}_{(\mathbf{\widetilde{A}}^{[n]})^T}\Big)\cb^{[n+1]}.
\end{align}

\subsection{Prediction-correction scheme}
For the construction of non-linear multiresolution analyses, we restrict ourselves to a special case of biorthogonal wavelet systems, namely prediction-correction schemes.
These schemes are typically associated with an interpolatory subdivision operator $\mathcal{S}_{\Ab^{[n]}}$ (predictor), i.e.\ an operator satisfying
$$\mathcal{D}_{\delta}\mathcal{S}_{\Ab^{[n]}}=\Db.$$

To obtain the other operators, we use the prediction-correction scheme as defined in \cite[Eq. (25)]{cotronei17} in terms of symbols:
\begin{equation}\label{eq:pred_correc_symbols}
 %\Ab^{[n]} \text{ interpolatory}, \quad 
 \Bb^{[n]}(z)=zI, \quad
 \mathbf{\widetilde{A}}^{[n]}(z)=\Db^{-1}, \quad
 \mathbf{\widetilde{B}}^{[n]}(z)=z\Db^{-1}(\Ab^{[n]})^{\sharp}(-z),
\end{equation}
again with notation $(\Ab^{[n]})^{\sharp}(z):=(\Ab^{[n]})^T(z^{-1})$.

%%%%%%%%%%%%%%%%%%%%%%%%%%%%%%%%%%%%%%%%%%%
% Prediction-correction in terms of operators
%%%%%%%%%%%%%%%%%%%%%%%%%%%%%%%%%%%%%%%%%
\begin{Lemma}\label{lem:pred_correc_operator}
The prediction-correction scheme defined in \eqref{eq:pred_correc_symbols} can be written in terms of operators in the following way:
\begin{enumerate}
%\item $\Ab^{[n]}$ interpolatory $\implies$ $\mathcal{D}_{\delta}\mathcal{S}_{\Ab^{[n]}}=D,$ where $\delta_i=\begin{cases}I & \text{if } i=0\\ 0 & \text{else} \end{cases}.$
\item
$\mathcal{S}_{\Bb^{[n]}}=\mathcal{L}^{-1}\mathcal{S}_{\delta},$ 
\item $\mathcal{D}_{(\mathbf{\widetilde{A}}^{[n]})^T}=\Db^{-1}\mathcal{D}_{\delta}$.
\item\label{tildeB}
%\begin{equation*}
$
 \mathcal{D}_{(\mathbf{\widetilde{B}}^{[n]})^T}=
\Db\,\mathcal{D}_{(\mathbf{\widetilde{A}}^{[n]})^T}\mathcal{L}\,
\Big(\operatorname{id}-\mathcal{S}_{\Ab^{[n]}}\mathcal{D}_{(\mathbf{\widetilde{A}}^{[n]})^T}\Big)=
\mathcal{D}_{\delta}\mathcal{L} \Big(\operatorname{id}-\mathcal{S}_{\Ab^{[n]}}\Db^{-1}\mathcal{D}_{\delta} \Big) 
$
%=\mathcal{D}_{\delta}\mathcal{L} \Big(\operatorname{id}-\mathcal{S}_{\Ab^{[n]}}\mathcal{D}_{(\mathbf{\widetilde{A}}^{[n]})^{T}} \Big)
%\end{equation*}
%This applies to the decomposition used in Definition \ref{def:decrec}.
% \item If we use the decomposition of \eqref{dec2}, then
% $$\mathcal{D}_{(\mathbf{\widetilde{B}}^{[n]})^T}=D\,\mathcal{D}_{(\mathbf{\widetilde{A}}^{[n]})^T}\mathcal{L}=\mathcal{D}_{\delta}\mathcal{L}.$$
\end{enumerate}
\end{Lemma}
%%%%%%%%%%%%%%%%%%%%%%%%%%%%%%%%%%%%%%%%%%%%%%
% REMARK
%%%%%%%%%%%%%%%%%%%%%%%%%%%%%%%%%%%%%%%%%%%%%%%
\begin{Remark}
From \Cref{lem:pred_correc_operator} it is apparent that the prediction-correction construction of \cite{cotronei17} is a Hermite version of \cite[Example 1.1.]{grohs12}.
\end{Remark}
%%%%%%%%%%%%%%%%%%%%%%%%%%%%%%%%%%%%%%%%%%%%%%%
% PROOF
%%%%%%%%%%%%%%%%%%%%%%%%%%%%%%%%%%%%%%%%%%%%%%%
\begin{proof}[Proof of \Cref{lem:pred_correc_operator}]
The first two parts are immediate from the definition of the symbol.

To see part (\ref{tildeB}), we first compute $(\mathbf{\widetilde{B}}^{[n]})^T_k$ from its symbol:
\begin{equation}
(\mathbf{\widetilde{B}}^{[n]})^T_k=(-1)^{1-k}\Ab^{[n]}_{1-k}\Db^{-1},
\end{equation}
see also \cite[p.~14]{cotronei17}. Therefore
\begin{equation}\label{op4}
(\mathcal{D}_{(\mathbf{\widetilde{B}}^{[n]})^T}\cb)_i=\sum_k(\mathbf{\widetilde{B}}^{[n]})^T_k\cb_{k+2i}=
\sum_k(-1)^{1-k}\Ab^{[n]}_{1-k}\Db^{-1}\cb_{k+2i}.
\end{equation}
Now compute the other operator, using the first two parts of this lemma and the interpolation property of $\Ab^{[n]}$:
\begin{align*}
&\Big(\Db\,\mathcal{D}_{(\mathbf{\widetilde{A}}^{[n]})^T}\mathcal{L}
\Big(\operatorname{id}-\mathcal{S}_{\Ab^{[n]}}\mathcal{D}_{(\mathbf{\widetilde{A}}^{[n]})^T}\Big)\cb\Big)_i=
(\mathcal{D}_{\delta}\mathcal{L}\cb)_i-(\mathcal{D}_{\delta}\mathcal{L}\mathcal{S}_{\Ab^{[n]}} \Db^{-1}\mathcal{D}_{\delta}\cb)_i\\
&=\cb_{2i+1}-(\mathcal{S}_{\Ab^{[n]}} \Db^{-1}\mathcal{D}_{\delta}\cb)_{2i+1} \\
&=\cb_{2i+1}-\sum_j \Ab^{[n]}_{-2j+1}\Db^{-1}(\mathcal{D}_{\delta}\cb)_{i+j}\\
&=\cb_{2i+1}-\sum_j \Ab^{[n]}_{-2j+1}\Db^{-1}\cb_{2(i+j)}\\
&=
%\overset{\text{interp.}}{=}
\cb_{2i+1}-\sum_j \Ab^{[n]}_{-2j+1}\Db^{-1}\cb_{2(i+j)}+\sum_j\Ab^{[n]}_{-2j}\Db^{-1}\cb_{2(i+j)+1}-\cb_{2i+1}\\
&=\sum_k(-1)^k\Ab^{[n]}_k\Db^{-1}\cb_{2i+1-k}\\
&=\sum_k(-1)^{1-k}\Ab^{[n]}_{1-k}\Db^{-1}\cb_{2i+k}.
\end{align*}
From \eqref{op4} the result follows.
\end{proof}

%%%%
% Decomposition and reconstruction scheme
%%%%
Based on \Cref{lem:pred_correc_operator}, the decomposition scheme (\Cref{def:dec}) in the prediction-correction case is given by
\begin{align}\label{pred-corr-decom}
\cb^{[n]}_i & = \Db^{-1}\cb^{[n+1]}_{2i} \\\nonumber
\db^{[n]}_i &= \left(\cb^{[n+1]}-\mathcal{S}_{\Ab^{[n]}}\cb^{[n]}\right)_{2i+1}.
\end{align}
We note that due to the interpolation property of $\mathcal{S}_{\Ab^{[n]}}$, i.e., since $\mathcal{S}_{\Ab^{[n]}}\cb^{[n]}_{2i}=\Db\cb^{[n]}_i=\cb^{[n+1]}_{2i}$, we have
$\left(\cb^{[n+1]}-\mathcal{S}_{\Ab^{[n]}}\cb^{[n]}\right)_{2i}=0$.

Given a function $f\in C^1(\RR,V)$, the discrete data $\cb^{[n]}_i$ is interpreted as samples of the function and its derivative at $i/2^n$. This means
\begin{equation}\label{eq:cn_function}
    \cb^{[n]}=\Db^n\vb_f^{[n]}.
\end{equation}
Through this interpretation, we obtain the \emph{Hermite wavelet transform} of $f$, which represents $f$ in terms of the decomposition sequence
\begin{equation}\label{eq:decom_seq}
    \cb^{[0]},\db^{[0]},\db^{[1]},\ldots
\end{equation}

The reconstruction scheme in the prediction correction case is
\begin{align}\label{pred-corr-recon}
 \cb^{[n+1]}_{2i} &= \Db\cb^{[n]}_i \\ \nonumber
 \cb^{[n+1]}_{2i+1} &= \left(\mathcal{S}_{\Ab^{[n]}}\cb^{[n]}\right)_{2i+1}+\db^{[n]}_i.
\end{align}
This can be used to reconstruction the function $f$ from the decomposition sequence \eqref{eq:decom_seq}.

%%%%%%%%%%%%%%%%%%%%%%%%%%%%%%%%%%%%%%%%%%%%%%%%%%%%%%%%%%%%%%%%%%
\section{Hermite subdivision and wavelets for manifold-valued data}\label{sec:manifoldSubDiv}
%%%%%%%%%%%%%%%%%%%%%%%%%%%%%%%%%%%%%%%%%%%%%%%%%%%%%%%%%%%%%%%%%%
%say something about how Hermite data looks like on manifolds, i.e. $TM$.

%%%%%%%%%%%%%%%%%%%%%%%%%%%%%%%%%%%%%%%%%%%%%%%%%
\subsection{Basic constructions in manifolds}\label{sec:manifold}
%%%%%%%%%%%%%%%%%%%%%%%%%%%%%%%%%%%%%%%%%%%%%%%%%
%
% which manifolds
By $M$ we denote a smooth, finite-dimensional manifold which carries a linear connection\footnote{By this we mean a linear connection on the tangent bundle $TM\to M$, which induces a covariant derivative in the sense of \cite[Section 19.11-19.12]{michor08}.}.
A linear connection allows to compute derivatives along tangent directions of vector fields (and more general, tensors), see \cite[Chapter IV]{michor08} for an introduction. The most important examples of such manifolds are Riemannian manifolds with the Levi-Civita connection \cite{doCarmo92}, and Lie groups with a Cartan-Schouten connection \cite{cogliati18,postnikov01}.

As $M$ carries a linear connection we have notions of parallel transport, geodesics and the exponential map, which we now define.

%% derivative along curve
By $TM$ we denote the tangent bundle, and by $T_pM$ the tangent space at $p\in M$, which is a linear space. For $I=[0,1]$, let $c:I \to M$ be a smooth curve such that $c(0)=p$ and $c(1)=q$ with $p,q \in M$. A \emph{vector field along $c$} is a smooth curve $\vfield:I \to TM$ such
that $\vfield(t)\in T_{c(t)}M$. Via the linear connection on $M$ we can differentiate vector fields along $c$. If in local charts (using Einstein summation) we have $\vfield=v^k\partial_k$ and $\dot{c}=x^k\partial_k$, then
\begin{equation*}
    \frac{D\vfield}{dt}:=\left(\frac{dv^k}{dt}+v^ix^j\Gamma_{ji}^k\right)\partial_k.
\end{equation*}
Here the coefficients $\Gamma_{ji}^k$ are uniquely determined by the underlying linear connection. If $M$ is a Riemannian manifold, they are called \emph{Christoffel symbols}.

The vector field $\vfield$ is called \emph{parallel along $c$} if
\begin{equation*}
    \frac{D\vfield}{dt}=0.
\end{equation*}
In charts this is a linear ODE, which implies that for a curve $c$, $c(0)=p$ and $v\in T_pM$ there exists a unique vector field $\vfield$ along $c$ such that $\vfield(0)=v$.

%% geodesic
Since $\dot{c}$ is a vector field along $c$, we define a \emph{geodesic} to be a curve $c$ satisfying
\begin{equation*}
    \frac{D\dot{c}}{dt}=0.
\end{equation*}
There exists a unique geodesic joining two points $p$ and $q$ (if not too far apart). In the Riemannian case, geodesics locally minimize length.

%% exponential map
The \emph{exponential map} is defined by $\exp_p(v):=g(1)$, where $g$ is the unique geodesic $g$ satisfying $g(0)=p$ and $\dot{g}(0)=v$.

%% dense enough input
We mention that the exponential map is always smooth, but in general not globally defined. Two important examples for which it is globally defined are complete Riemannian manifolds and matrix groups \cite{helgason79,onishchik93}. Similarly, the inverse exponential is generally only smooth if $p$ and $q$ are close together. Manifold-valued subdivision schemes often rely on the exponential map and therefore results are usually only valid for ``dense enough'' input data, see for example \cite{wallner05,wallner11,grohs10,moosmueller16,xie07}. However, there exist convergence results valid for all input data in specific cases \cite{huening20,huening19,wallner11}. Dense enough input data is also a necessary assumption for our results in \Cref{sec:results}.

%% parallel transport
If $c(0)=p$ and $c(1)=q$, then the \emph{parallel transport along $c$} is the linear map $P_p^q(c): T_pM \to T_qM$, $v \mapsto \vfield(1)$, where $V$ is the unique parallel vector field along $c$ with $\vfield(0)=v$. The map $P_p^q(c)$ is an isomorphism, and if $M$ is a Riemannian manifold, it is also an isometry.
The parallel transport satisfies
\begin{equation}\label{eq:parallel}
    P_m^q(c) \circ P_p^m(c) =P_p^q(c),
\end{equation}
where $m$ is a point on $c$.
In this paper we always choose the curve to be the geodesic joining $p$ and $q$ when we compute the parallel transport. We introduce the simplified notation 
\begin{equation*}
[v]_q:=P_p^q(g)(v),
\end{equation*}
where $v\in T_pM$ and $g$ is the geodesic from $p$ to $q$.
\Cref{eq:parallel} now reads $[[v]_m]_q=[v]_q$.

%%%%%%%%%%%%%%%%%%%%%%%%%%%%%%%%%%%%%%%%%%%%%%%%%%%%%%%%%%%%%
\subsection{Hermite subdivision schemes for manifold-valued data and the proximity condition}
%%%%%%%%%%%%%%%%%%%%%%%%%%%%%%%%%%%%%%%%%%%%%%%%%%%%%%%%%%%%%
%%% INTRO

%%%%% Subdivision operator on M
Following \cite{moosmueller16}, we define a Hermite subdivision operator for manifold-valued data.
%% definition of manifold-valued SD operator
\begin{Definition}\label{def:Hermite_sd_manifold}
A \emph{Hermite subdivision operator on $M$} is a map $\cT: \ell(\ZZ,T M) \to \ell(\ZZ, T M)$ such that
\begin{enumerate}
    \item $\cL^2 \cT = \cT \cL$, where $\cL$ is the left shift operator 
    \eqref{eq:shift},
    \item $\cT$ has compact support, i.e.\ there exists $N$ such that $(\cT \cb)_{2j}$ and $(\cT \cb)_{2j+1}$ depend only on $\cb_{j-N},\ldots,\cb_{j+N}$, for all $j \in \ZZ$ and $\cb\in \ell(\ZZ,T M)$.
\end{enumerate}
\end{Definition}
Compare this definition with the properties of linear Hermite subdivision operators \eqref{eq:subdivisionOperator} and \eqref{eq:sub-wavelet-shift-commute}.

% Based on \Cref{def:Hermite_sd_manifold}, a Hermite subdivision scheme on $M$ is 

%%% Subdivision operator M from linear operator
We use a linear Hermite subdivision operator $\cS_{\Ab}$, with mask $\Ab$ of the form \eqref{eq:typeMasks}, to define a manifold-valued analogue $\cT_{\Ab}$ satisfying the properties of \Cref{def:Hermite_sd_manifold}. This is based on the parallel transport construction of \cite{moosmueller17}.

Choose a base point sequence $m\in \ell(\ZZ,M)$. For $\cb = (p, v)^T \in \ell(\ZZ,T M)$ we define
\begin{equation}\label{eq:manifold_so_from_linear}
    \left(\cT_{\Ab}\cb\right) = \tilde{\bf{c}},
\end{equation}
where $\tilde{\cb} = (\tilde{p},\tilde{v})^T\in \ell(\ZZ,T M)$ is given by
\begin{align*}
    \tilde{p}_j &= \exp_{m_j}
    \left( \sum_{k\in \ZZ} a_{j-2k}^{00} \exp^{-1}_{m_j}(p_k)
    + a_{j-2k}^{01} [v_{k}]_{m_j}
    \right), \\
    \tilde{v}_j &= \left[\sum_{k\in \ZZ} a_{j-2k}^{10} \exp^{-1}_{m_j}(p_k)
    + a_{j-2k}^{11} [v_{k}]_{m_j}
    \right]_{\tilde{p}_j}
\end{align*}
for $j \in \ZZ$.

%%%% Subdivision scheme on manifold
From the manifold-valued subdivision operator based on a mask $\Ab$ \eqref{eq:manifold_so_from_linear}, we can define a manifold-valued subdivision scheme as the iterative process to construct $\cb^{[n]} \in \ell(\ZZ,TM)$ from $\cb^{[0]}\in \ell(\ZZ,TM)$ via
\begin{equation}
\label{eq:subdivisionScheme_manifold}
\Db^{\itVar+1}\cb^{[\itVar+1]} = \mathcal{T}_{\Ab^{[\itVar]}}\Db^{\itVar} \cb^{[\itVar]},
\quad \itVar \in \NN,
\end{equation}
where $(\Ab^{[n]},n\in \NN)$ is a sequence of masks.

%%%% Charts

Results for manifold-valued subdivision schemes on topics such as convergence, smoothness, and approximation order, are often derived from their linear counterparts via a proximity condition \cite{moosmueller16,moosmueller17,wallner11, xie07,grohs10,grohs07,wallner05}. A comparison between a linear and a manifold-valued operator only makes sense in a chart or an embedding of $M$. In this paper we use charts and thus assume that $TM\subset V^2$. 

We now define a proximity condition for Hermite subdivision operators as in \cite{moosmueller16}, which is also to be understood in charts.

%%%%% Proximity condition
\begin{Definition}[Proximity condition]\label{def:proximity}
Let $\left(\cS_{\Ab^{[n]}}: n  \in \NN \right)$ be a sequence of linear Hermite subdivision operators. Let $\left(\cT_{\Ab^{[n]}}: n  \in \NN \right)$ be its manifold-valued analogue defined via \eqref{eq:manifold_so_from_linear}. The proximity condition is satisfied if there exists a constant $C$ such that
\begin{equation*}
    \left\|\left(\cS_{\Ab^{[n]}}- \cT_{\Ab^{[n]}}\right) \begin{pmatrix} p \\ v \end{pmatrix} \right\|_{\infty}
    \leq C \,\left\|\begin{pmatrix} \Delta p \\ v \end{pmatrix}\right\|_{\infty}^2,
    %\|\cG^{[n]} \cb\|^2,
    \quad n\in \NN, (p,v)^T \in \ell(\ZZ,TM),
\end{equation*}
\end{Definition}

%%% Results on proximity and convergence

In \cite[Corollary 1]{moosmueller17} it is shown that if the base point sequence is chosen as either $m_i=p_i$ or as the geodesic midpoint between $p_i$ and $p_{i+1}$, and the input data is bounded, then the proximity condition between $\cS_{\Ab}$ and $\cT_{\Ab}$ is satisfied. Therefore, in this paper, we choose the base point sequence as either one of those sequences.

\subsection{Manifold-valued prediction-correction scheme}\label{sec:manifold-prediction-correction}
We define operations $\oplus$ and $\ominus$ in manifolds as generalization of $+,-$ in vector spaces. Indeed, the operations we define are extensions of $\oplus,\ominus$ defined in \cite{wallner11,grohs12} for point-data to Hermite data. 

%The operation $\oplus:TM\times (TM\oplus TM) \rightarrow TM$ is given by
We consider point-vector Hermite data $(p,v)^T$ and vector-vector data $(u_0,u_1)^T$, which is an element of $T_qM\oplus T_qM$, with $q\in M$, hence an element of a fiber of $TM\oplus TM$. We define the addition of such elements as:
\begin{align}\label{eq:f-plus}
%&\oplus:TM\times (TM\oplus TM) \rightarrow TM\\\nonumber
    {p\choose v}&\oplus {u_0 \choose u_1} := %\\ \nonumber
    {\exp_p([u_0]_p) \choose [v]_{\exp_p([u_0]_p)}+[u_1]_{\exp_p([u_0]_p)}}.
\end{align}
Similarly, for point-vector data $(p,v)^T,(q,u)^T$ we define their difference as
%while the operation $\ominus:TM\times TM\rightarrow TM\oplus TM$ is defined as
\begin{align}\label{eq:g-minus}
%    &\ominus:TM\times TM\rightarrow TM\oplus TM \\\nonumber
    {q \choose u}& \ominus {p \choose v}  := {\exp_p^{-1}(q) \choose [u]_p-v }.
\end{align}
The resulting element lies in the fiber $T_pM\oplus T_pM$.
In \Cref{lem:plus_minus_ok} below we show that these operations satisfy similar properties as the operations on point-data defined in \cite{wallner11}.

%%%%%%%%%%%%%%%%%%%%%%%%%%%%%%%%%%%%%%%%%%
% LEMMA: OPERATIONS MAKE SENSE
%%%%%%%%%%%%%%%%%%%%%%%%%%%%%%%%%%%%%%%%%%
\begin{Lemma}\label{lem:plus_minus_ok}
Consider point-vector data $\ab, \tilde{\ab}$ and vector-vector data $\bb$.
Then we have the following properties:
\begin{align*}
\ab\oplus(\tilde{\ab}\ominus \ab)&=\tilde{\ab},\\
(\ab\oplus \bb)\ominus \ab&=[\bb]_p,
\end{align*}
with $[\bb]_p = ([u_0]_p,[u_1]_p)^T$ when $\bb = (u_0,u_1)^T$.
\end{Lemma}
%%%%%%%%%%%%%%%%%%%%%%%%%%%%%%%%%%%%%%%%%%%
\begin{proof}
Let $\ab =(p,v)^T, \tilde{\ab} = (\tilde{p},\tilde{v})^T$ and $\bb =(u_0,u_1)^T$.
% 
% Furthermore introduce the short hand notation $v=(v_1,\ldots,v_d), \tilde{v}=(\tilde{v}_1,\ldots,\tilde{v}_d), u=(u_1,\ldots,u_d)$ and $[v]_p=([v_1]_p,\ldots,[v_d]_p)$, $p\in M$.
Then \eqref{eq:g-minus} implies
\begin{align*}
\tilde{\ab}\ominus \ab
= (\exp_p^{-1}(\tilde{p}),[\tilde{v}]_p-v)^T.
\end{align*}
From \eqref{eq:f-plus} we see that the first entry of $\ab\oplus(\tilde{\ab}\ominus \ab)$ is $\tilde{p}$ and
\begin{align*}
\ab\oplus(\tilde{\ab}\ominus \ab)
= (\tilde{p},[v]_{\tilde{p}} +[[\tilde{v}]_{p}-v]_{\tilde{p}})
=(\tilde{p},[\tilde{v}]_{\tilde{p}})=\tilde{\ab}.
\end{align*}
Similarly, \eqref{eq:f-plus} and \eqref{eq:g-minus}
% \begin{align*}
%     \cb\oplus \db = (p \oplus u_0, [v]_{p \oplus u_0}+ [u]_{p \oplus u_0}).
% \end{align*}
\begin{align*}
    (\ab\oplus \bb)\ominus \ab = ([u_0]_p, v+[u_1]_p-v)
    = ([u_0]_p,[u_1]_p)=[\bb]_p.
\end{align*}
This concludes the proof.
\end{proof}
\begin{Remark}
 If $\ab$ and $\bb$ are taken from the same fiber, i.e.\ $v,u_0,u_1 \in T_pM$, then $(\ab\oplus \bb)\ominus \ab=\ab$. 
\end{Remark}

Based on $\oplus, \ominus$ and \eqref{pred-corr-decom}, \eqref{pred-corr-recon}, we can define a prediction-correction scheme for manifold-valued Hermite data where the decomposition scheme is
\begin{align}\label{pred-corr-decom_mani}
\cb^{[n]}_i & = \Db^{-1}\cb^{[n+1]}_{2i} \\\nonumber
\db^{[n]}_i &= \left(\cb^{[n+1]}\ominus \mathcal{T}_{\Ab^{[n]}}\cb^{[n]}\right)_{2i+1}.
\end{align}
Similar to \eqref{eq:cn_function}, for a function $f\in C^1(\RR,M)$, we can interpret $\cb^{[n]}=\Db^{n}\vb_f^{[n]}$ and use \eqref{pred-corr-decom_mani} as the decomposition sequence of $f$.
The reconstruction scheme is then defined by 
\begin{align}\label{pred-corr-recon_mani}
 \cb^{[n+1]}_{2i} &= \Db\cb^{[n]}_i \\ \nonumber
 \cb^{[n+1]}_{2i+1} &= \left(\mathcal{T}_{\Ab^{[n]}}\cb^{[n]}\right)_{2i+1}\oplus\db^{[n]}_i.
\end{align}

%%%%%%%%%%%%%%%%%%%%%%%%%%
%%%%%%%%%%%%%%%%%%%%%%%%%
\subsection{Coefficient decay for manifold-valued Hermite wavelets}\label{sec:manifold-decay}
%%%%%%%%%%%%%%%%%%%%%%%%%

We now generalize the linear wavelet coefficient decay result of \cite{cotronei19} to the manifold-valued case.

%%%%%%%%%%%%%%%%%%%%%%%%%%%%%%%%%%
% THEOREM: COEFFICIENT DECAY
%%%%%%%%%%%%%%%%%%%%%%%%%%%%%%%%%%
\begin{Theorem}\label{manifold_decay}
	Let $S({\Ab}^{[n]}: n \geq 0)$ be a $C^1$-convergent interpolatory
	Hermite subdivision scheme satisfying the $W$-spectral
	condition.
	% and whose 
	% basic limit functions $(\Fb^{[n]}: n \geq 0)$ generate an interpolatory MRA. 
	Moreover assume that there exists $N \in \NN$ such that
	$\supp(\Ab^{[n]}) \subseteq [-N,N]$ for all $n \in \NN$,
	and that $\sup_{n \in \NN}\|\Fb^{[n]}\|_{\infty} < \infty$.
	Let $M$ be a manifold (as described in \Cref{sec:manifold}) and let $f \in C^{1}_u(\RR,M)$. We assume that $\cb^{[N]}$ is dense enough. Then the associated manifold-valued wavelet coefficients
	$\db^{[n]}$ \eqref{pred-corr-decom_mani} satisfy the following property:
	For $R<1$, there exist $m \in \NN$ and a constant $C>0$, depending on
	$W,R,f,N,M$ and the subdivision scheme, such that
	\begin{equation*}
	%\label{eq:leveldep_decay}
	\|\db^{[n]}\|_{\infty}\leq C\, 2^{-2n}, \qquad n \geq m.
	\end{equation*}
\end{Theorem}

%%%%%%%%%%%
% PROOF
%%%%%%%%%%%
\begin{proof}
We first note that for bounded sequences $\ab,\bb \in \ell_{\infty}(\ZZ,TM)$, the operator $\ominus$, as defined in \eqref{eq:g-minus}, satisfies
\begin{equation*}
    \left\|\ab \ominus \bb\right\|_{\infty} \leq C 
    \left\|\ab-\bb\right\|_{\infty},
\end{equation*}
for some constant $C$.
This follows from the linearizations $\exp^{-1}_p(q)=q-p+O(\|q-p\|^2)$ and $P_q^p(u) = u + O(\|q-p\|\|u\|)$ for $q\to p$ and fixed $u$, compare \cite[Lemma 1]{moosmueller17}.
Therefore, we have
\begin{align}\nonumber
\|\db^{[\itVar]}\|_{\infty} = & \| \cb^{[\itVar+1]}\ominus {\mathcal{T}}_{\Ab^{[\itVar]}}\cb^{[\itVar]}\|_{\infty}
\le C
\| \cb^{[\itVar+1]}- {\mathcal{T}}_{\Ab^{[\itVar]}}\cb^{[\itVar]}\|_{\infty}\\ \label{eq:d_bound}
&\le C\left(
\| \cb^{[\itVar+1]}- {\mathcal{S}}_{\Ab^{[\itVar]}}\cb^{[\itVar]}\|_{\infty}+
\| {\mathcal{S}}_{\Ab^{[\itVar]}}\cb^{[\itVar]}- {\mathcal{T}}_{\Ab^{[\itVar]}}\cb^{[\itVar]}\|_{\infty}
\right),
\end{align}
The first part is bounded by $C\,2^{-2\itVar}$ whenever $n\geq m$ by the linear wavelet decay result of \cite[Theorem 11]{cotronei19}.
%see also Theorem \ref{T:leveldep_decay} in the appendix.
For the second part, the proximity condition (\Cref{def:proximity}) implies:
\begin{align*}
\left\| {\mathcal{S}}_{\Ab^{[\itVar]}}\cb^{[\itVar]}- {\mathcal{T}}_{\Ab^{[\itVar]}}\cb^{[\itVar]}\right\|_{\infty}
\leq \left\|\left( \begin{array}{cc}
    \Delta & 0 \\
     0 &  1
\end{array} \right)
\cb^{[\itVar]}
\right\|_{\infty}^2.
\end{align*}
Since $\cb^{[\itVar]} = \Db^{\itVar} \vb_{f}^{[\itVar]}$, the two component of the right side are given by
\begin{equation*}
\left(
   \begin{array}{cc}
    \Delta & 0 \\
     0 &  1
\end{array}
 \right) \cb^{[\itVar]}=
\left(
\begin{array}{c}
     \Delta f(j/2^{\itVar})  \\
     2^{-\itVar} f'(j/2^{\itVar})
\end{array}
\right).
\end{equation*}
Since $f\in C^1_u(\RR,M)$, $f'$ is bounded and therefore $f$ is Lipschitz. Thus 
$$\left\|\Delta f\left(\frac{j}{2^{\itVar}}\right)\right\|_{\infty}=
\left\|f\left(\frac{j+1}{2^{\itVar}}\right)-f\left(\frac{j}{2^{\itVar}}\right)\right\|_{\infty} \leq C \,{2^{-\itVar}}$$
and we obtain the bound
\begin{equation*}
\left\| {\mathcal{S}}_{\Ab^{[\itVar]}}\cb^{[\itVar]}- {\mathcal{T}}_{\Ab^{[\itVar]}}\cb^{[\itVar]}\right\|_{\infty} \leq 
    \left\|\left( \begin{array}{cc}
    \Delta & 0 \\
     0 &  1
\end{array} \right)
\cb^{[\itVar]}
\right\|_{\infty}^2 \leq C\,2^{-2\itVar}.
\end{equation*}
This bound together with the estimate \eqref{eq:d_bound} and the linear wavelet coefficient result \cite[Theorem 11]{cotronei19} concludes the proof.
\end{proof}
%%%%%%%%%%%%%%%%%%%%%%%%%%%%%%%%%%%%%%%%%%%
% CONCLUSION
%%%%%%%%%%%%%%%%%%%%%%%%%%%%%%%%%%%%%%%%%%
%%%%%%%%%%%%%%%%%%%%%%%%%%%%%%%%%%%%%%%%%%
\section{Conclusions}
%%%%%%%%%%%%%%%%%%%%%%%%%%%%%%%%%%%%%%%%%

In this paper we have provided a framework for the construction of Hermite-type multiwavelets in a manifold setting. In particular we have extended to such a setting a recent result  about the decay of the wavelet coefficients  \cite{cotronei19}. Our ideas go in the direction of providing efficient representations of Hermite manifold-valued data  as in a  traditional wavelet analysis, for example for  compression or denoising applications. Future research will focus on such applications and on the generalization of the obtained theoretical results to the case of higher order derivatives.

%%%%%%%%%%%%%%%%%%%%%%%%%%%%%%%%%%%%%%%%%%

%%%%%%%%%%%%%%%%%%%%%%%%%%%%%%%
% FUNDING
%%%%%%%%%%%%%%%%%%%%%%%%%%%%%%
\section*{Acknowledgments}
Mariantonia Cotronei is member of  RITA (Research ITalian network on Approximation), and of  INdAM-GNCS and UMI-TAA research groups.
\par\smallskip\noindent
Nada Sissouno and Caroline Moosm{\"u}ller acknowledge partial funding from an Entrepreneurial Award in the Program ``Global Challenges for Women in Math Science" funded by the Faculty of Mathematics at the Technical University of Munich.
\par\smallskip\noindent
Caroline Moosm{\"u}ller is supported by NSF award DMS-2111322.

%%%%%%%%%%%%%%%%%%%%%%%%%%%%%%%%%%%%%%%%%%%%%%%%%%%%%%
% BIB
%%%%%%%%%%%%%%%%%%%%%%%%%%%%%%%%%%%%%%%%%%%%%%%%%%%%%%
% \bibliographystyle{abbrv}
% \bibliography{lit}
%%%%%%%%%%%%%%%%%%%%%%%%%%%%%%%%%%%%%%%%%%%%%%%%%%%%%%%
% add bibliography directly for arXiv upload

\end{document}